\documentclass[a4paper,12pt]{article}

\usepackage{amsmath,amsfonts,latexsym,theorem,
}
\usepackage{pdfpages}
\usepackage{graphicx}
\usepackage{subfigure}
\usepackage{caption}
\nonstopmode

\numberwithin{equation}{section}

\newtheorem{Theorem}{Theorem}[section]
\newtheorem{Definition}{Definition}[section]
\newtheorem{Proposition}{Proposition}[section]
\newtheorem{Lemma}{Lemma}[section]

\newtheorem{Corollary}{Corollary}[section]
\newenvironment{Proofc}[1]{\smallskip\par\noindent\textsc{#1}\quad}%
  {\hfill$\Box$\bigskip\par}
\newenvironment{Proof}{\begin{Proofc}{Proof}}{\end{Proofc}}

\newtheorem{Remark}{Remark}[section]

\def\d{\delta}

\def\g{\gamma}
\def\G{\Gamma}
\def\l{\lambda}

\def\og{\omega}
\def\O{\Omega}

\def\e{\varepsilon}

\def\pd{\partial}

\def\half{\frac{1}{2}}

\newcommand{\cA}{{\cal A}}

\newcommand{\cX}{{\cal X}}

\newcommand{\R}{{\mathbb R}}
\newcommand{\N}{{\mathbb N}}

\def\pd{\partial}

\begin{document}
\title{ Eikonal equations on the Sierpinski gasket}

\author{Fabio Camilli\footnotemark[1] \and Raffaela Capitanelli\footnotemark[1] \and Claudio Marchi\footnotemark[2] }

\date{version: \today}
\maketitle

\footnotetext[1]{Dip. di Scienze di Base e Applicate per l'Ingegneria,  ``Sapienza'' Universit{\`a}  di Roma, via Scarpa 16,
 00161 Roma, Italy, ({\tt e-mail:camilli,capitanelli@dmmm.uniroma1.it})}
\footnotetext[2]{Dip. di Matematica, Universit\`a di Padova, via Trieste 63, 35121 Padova, Italy ({\tt marchi@math.unipd.it}).}

\begin{abstract}
We study the eikonal equation on the Sierpinski gasket in the spirit of the construction of the Laplacian in Kigami \cite{k1}: we consider   graph eikonal equations on the prefractals  and we show that the solutions of these problems converge  to a function  defined on the fractal set. We characterize this limit function as the unique metric viscosity solution to the eikonal equation on the Sierpinski gasket according to the definition introduced in \cite{ghn}.
\end{abstract}
 \begin{description}
\item [{\bf MSC 2000}:]  35R02, 49L25, 28A80.
\item [{\bf Keywords}:] Sierpinski gasket,  eikonal equation, viscosity solution.
\end{description}
\section{Introduction}

In this paper we consider the eikonal equation
\begin{equation}\label{HJi}
    |Du|=f(x) \qquad \textrm{in }  S,
\end{equation}
where $S$ is the Sierpinski gasket. The eikonal equation is the prototype for  a general  class of first Hamilton-Jacobi equations with convex Hamiltonian (see \cite{bcd}) and, moreover, it arises in several applications in connection with geometric optics, wave front propagation, interfaces evolution, granular matter theory, etc.\par
The analysis of differential operators on irregular sets such as the Sierpinski gasket has been developed
since the late 80's, starting with the pioneering works of Goldstein \cite{g}, Kusuoka \cite{ku} and Lindstr{\o}m \cite{l}. Since there is no natural notion of derivative on general closed sets, the notion of Laplacian on fractals  is introduced by approximating the set from within and performing a limiting process.
The probabilistic version of this approach was introduced by Kusuoka \cite{ku} and Lindstr{\o}m \cite{l} who  considered suitable scaled random walks  on the prefractal  and then passed to the limit as the graph approaches the fractal so to define a   Brownian motion on $S$. The corresponding  analytical approach  was    taken by Kigami (see \cite{k1} and references therein)
who defined   $\Delta u$ on the class of post-critically  finite  (pcf in short) fractals as the uniform  limit of  suitable scaled finite difference schemes  on the prefractal.
On pcf fractals the two approaches give rise to the same self-adjoint differential operator
which is identified as the Laplacian  (see \cite{k1}, \cite{s2}). The interior approximation preserves  some of the classical properties satisfied by the Laplacian in the  Euclidean setting, but other properties are typical of the fractal one (\cite{m}). \par
In this paper we     pursue  an interior approximation   approach on the Sierpinski gasket $S$ for the eikonal equation \eqref{HJi}.
The starting point of our analysis is the work of Manfredi, Oberman and Sviridov \cite{mos} who consider  a general class of nonlinear elliptic difference equations on graphs.
This class includes the graph Laplacian, which coincides  up to a rescaling factor with the operator considered in \cite{k1}, and the graph eikonal equation.
The principle underlying the definition of Laplacian and harmonic functions on the Sierpinski gasket is the minimization of the energy. For the eikonal equation the  principle  is the optimal control interpretation of the solution at all the different levels (discrete and continuous on the prefractal and continuous on the fractal).   \par We consider the graph eikonal equation  on the prefractal $S^n$ with a Dirichlet boundary condition on the vertices of the simplex generating the Sierpinski gasket;  we characterize the unique solution of  the problem by a representation formula which,  for  $f\equiv 1$ and null boundary data, gives the  minimal vertex distance from the boundary. Moreover we infer that the sequence of the solutions of the graph eikonal equations on $S^n$ is compact hence, up to a subsequences, converges uniformly to a continuous function $u$ defined on $S$. By classical stability properties of viscosity solution theory it is only reasonable to expect the function $u$ to be the solution in some appropriate sense of the eikonal equation \eqref{HJi} on $S$.\par
In \cite{ghn}, the authors give a definition of viscosity solution for the eikonal equation \eqref{HJi} in a  general metric space $\cX$, which is consistent with the usual notion in the Euclidean space. The unique solution of the associated Dirichlet problem is characterized by a representation formula of control type.
This definition obviously applies to the Sierpinski gasket endowed with the metric induced by the path distance. We   show that the limit
of the sequence of the solutions of the graph eikonal equations on the prefractals converges to the solution defined as in \cite{ghn} on  $S$. To establish this link, we introduce a class of continuous eikonal equation on the prefractal $S^n$ and we obtain  an uniform estimate of the distance between the solutions of the discrete and the continuous eikonal equations on $S^n$. Then, thanks to the stability property of the definition in  \cite{ghn}, we are able to pass to the limit and to show that the solution  of the continuous eikonal equation  (and therefore also of the  graph eikonal equation)  on $S^n$ converges  to the solution of \eqref{HJi} on $S$. Note that  this intermediate step can be also seen as a sort of homogenization result for the continuous eikonal equation on the Sierpinski gasket.\par
The crucial point in \cite{ghn} is  a notion of metric derivative $|\xi'(t)|$   for a given  path $\xi=\xi(t)$ in $\cX$ although in general $\xi'(t)$ may be  not well defined. It follows that this elegant  theory is    confined to the case of the eikonal equation and  difficult to extend to a more general class of  Hamilton-Jacobi equation on $S$. Moreover the notion of viscosity supersolution in \cite{ghn} is not local, even if consistent with the Euclidean viscosity solution.
Instead we are able to extend  the interior approximation   approach to a  more general  class of Hamiltonians showing  that the corresponding sequence of the solutions of the graph Hamilton-Jacobi equations converges uniformly to a function $u$ on $S$. Lacking a definition of viscosity solution for general convex Hamiltonian on the Sierpinski gasket, our approach can be seen as a constructive way to define the solution to Hamilton-Jacobi equations on $S$. From this point of view  the previous result for the eikonal equation can be also  interpreted as a test  that the construction gives the correct solution on the limit fractal.
\par
The paper is organized as follows.
In Section \ref{sect:SG} we recall the definition of Sierpinski gasket. Sections \ref{s3} and  \ref{s4} are devoted to the study of the discrete
and respectively  continuous eikonal equation on the prefractal. In Section \ref{HJSierpinski} we study the eikonal equation on the Sierpinski gasket and prove the
convergence of the problems on the prefractal. Finally, in Section \ref{conclusion}  we consider some possible extensions of the results here developed.

\section{The Sierpinski gasket}\label{sect:SG}

Consider a unit regular $D$-dimensional simplex of vertices $\G=\{a_1,\dots,a_{D+1}\}$ in $\R^D$ (e.g., for  $D=1,2,3$, the simplex is respectively an interval, an equilateral triangle and a regular tetrahedron), and the $D+1$ mappings $\psi_i:\R^D\to\R^D$ defined by
\[\psi_i(x):=a_i+\half(x-a_i)\qquad i=1,\dots, D+1.\]
Iterating the $\psi_i$'s, we get the set
\[\displaystyle{\G^\infty=\cup_{n=0}^\infty \G^n}\]
where $\G^0\equiv \G$ and each $\G^n$ is given  by the union of the images of    $\G$ under the action of the maps  $\psi_{i_1}\circ\dots\circ\psi_{i_n}$
with $i_h\in \{1,\dots, D+1\}$, $h=1,\dots,n$. Then the Sierpinski gasket  $S$ is the closure of $\G^\infty$ (with respect to the Euclidean topology) and it is the unique non empty compact set $K$ which satisfies
\[K=\cup_{i=1}^{D+1}\psi_i(K).\]

For any $n$, we can identify $\G^n$ with the graph $(V^n,\sim_n)$, where $V^n=\G^n$ and $\sim_n$ is the following relation on $V^n$: for $x,y\in V^n$, $x\sim_n y$ if and only if the segment connecting $x$ and $y$ is the image of a side of the starting simplex  under the action of some $\psi_{i_1}\circ\dots\circ\psi_{i_n}$. When there is no possibility of confusion, for simplicity we shall write ``$\sim$'' instead of ``$\sim_n$''. Moreover, the graph $(V^n,\sim)$ can be naturally embedded in the network $S^n:=V^n\cup E^n$, where $E^n$ is the set of all segments of endpoints $x$, $y$ with $x,y\in V^n$ and $x\sim_n y$. In other words, $S^n$ is formed by the union of the images of the initial simplex under all maps $\psi_{i_1}\circ\dots\circ\psi_{i_n}$.
Hence $S=V\cup E$ where $V:=\cup_{i=1}^\infty V^n$ and $E:=\cup_{i=1}^\infty E^n$.
The $D+1$ vertices of the initial simplex assume the role of the boundary of $S^n$ and of the curve $S$ itself:
\[\partial S:=\G;\]
on this set a boundary condition will be prescribed.\par
We note that, for any $x\in S\setminus \G$, there are only two possibilities:
\begin{itemize}
  \item either  $x\in V\setminus \G$, hence, there exists $n\in \N$   and $v\in V^n\setminus \G$ such that $x=v$; in particular, $x$ is a ramification point with $4$ incident  edges in $S^m$, for any $m\ge n$;
  \item or $x \in   S\setminus (V\cup \G)$, hence  for $n$ sufficiently large, there exists $e_n\in E^n$ such that $x\in e_n$ and $x$ is not a ramification point in any $S^n$ (however, it is the limit of ramification points).
\end{itemize}

\begin{figure}
\centering
\subfigure[$\G$]{\includegraphics[scale=.4]{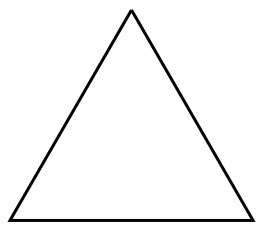}}
\subfigure[$\G^1$]{\includegraphics[scale=.4]{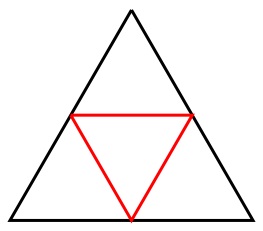}}
\subfigure[$\G^2$]{\includegraphics[scale=.4]{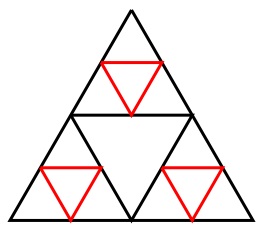}}
\subfigure[$\G^4$]{\includegraphics[scale=.3]{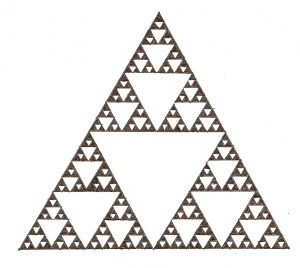}}
\end{figure}
We define the geodesic distance on  the network $S$.
\begin{Definition}\label{dist}
 For any $x,y\in S$,   set
 \begin{equation}\label{geodist}
   d(x,y):=\inf\left\{\ell(\g):\,\text{$\g$ is a path joining $x$ to $y$ }\right\}
 \end{equation}
 where $\ell(\g)$ is the length of $\g$.
\end{Definition}
We refer the reader to \cite{hs,s1} (see also \cite{k2}) for a complete characterization of the geodesic distance on the Sierpinski gasket. See \cite{hs} for the next result
\begin{Proposition}
For any $x,y\in S$, there holds
        \[
            |x-y|\le d(x,y)\le 1.
        \]
Moreover $d$ is continuous and the $\inf$ in \eqref{geodist} is a minimum.
\end{Proposition}

Since now on, we shall require the following hypotheses
\vskip 2mm
\noindent{\bf Assumptions}: $D=2$.
The function $f:S\to\R$ is continuous with
\[
\l:=\inf_Sf(x)>0.
\]

\section{A discrete eikonal  equation on the prefractal}\label{s3}

In this section, for $n\in\N$ fixed, we tackle the discrete eikonal equation on the graph $(V^n,\sim)$. We denote $V^n_0:=V^n\setminus \G$ the set of {\it interior} vertices while $\G$ will stand for the set of boundary vertices.
\par
 The distance between two adjacent vertices $x,y$ is
\[d_n(x,y)=\frac{1}{2^n}=:h_n\]
 while the distance for two arbitrary vertices  $x,y\in V^n$ is the {\it vertex distance}
\begin{equation}\label{geodistD}
d_n(x,y):=\inf\{d(x_0,x_1)+d(x_1,x_2)+\dots+d(x_{N-1},x_N) \}
\end{equation}
where the infimum is taken over all the finite path $\{x_0=x,x_1,\dots, x_N=y\}$ with $x_i\sim x_{i+1}$, $i=0,\dots,N-1$ connecting $x$ to $y$. In fact, since $V^n$ is a finite set, the ``$\inf$'' is a ``$\min$''.
Note that $d_n$ is the restriction  of the geodesic distance $d$ defined in \eqref{geodist} to $V^n$.
\par
For a function $u:V^n\to\R$,  we consider the discrete eikonal equation   (see \cite{mos})
\begin{equation}\label{HJD}
      |Du(x)|_n=f(x)   \qquad \textrm{in }  V_0^n
\end{equation}
where
\[
|Du(x) |_n:= \max_{y \sim x}\left\{-\frac{1}{h_n} (u(y)-u(x))\right\}.
\]

\begin{Definition}
 A  function $u:V^n\to\R$ is said a subsolution (respectively, a supersolution) of  \eqref{HJD} if
 \begin{equation*}
 |Du(x)|_n\le f(x)\quad (\text{resp.}, \, |Du(x)|_n\ge f(x))\qquad \forall x\in V_0^n.
 \end{equation*}
 A function $u$ is said a solution of   \eqref{HJD} if it is both a sub- and a supersolution of \eqref{HJD}.
 \end{Definition}
In the next proposition we prove  a  comparison theorem for \eqref{HJD}.
 \begin{Proposition}\label{comparisonD}
Let $u$ and $v$ be  a subsolution and, respectively,  a supersolution of \eqref{HJD}  such that $u\le v$ on $\G$. Then $u\le v$ in $V^n$.
 \end{Proposition}
 \begin{Proof}
We shall proceed using some ideas of \cite{mos}. Assume by contradiction
 that $\max_{V^n}\{u-v\}=\d>0$ and set $W=\arg\max_{V^n}\{u-v\}$, $m=\min\{v(x):\, x\in W\}$. \par
Let  $x\in W$ be  such that $v(x)=m$; in particular, $x$ belongs to $V_0^n$. Let $z\in V^n$ be such that $z\sim x$ and $|Dv(x)|_n=-h_n^{-1}(v(z)-v(x))$. By the  definition of sub- and supersolution we have
\begin{eqnarray*}
-\frac{1}{h_n}(v(z)-v(x))-f(x)&=&|Dv(x)|_n-f(x)\ge 0\ge |Du(x)|_n-f(x)\\
&\ge& -\frac{1}{h_n}(u(z)-u(x))-f(x);
\end{eqnarray*}
hence  $u(z)-v(z)\ge u(x)-v(x)=\d$ and therefore $z\in W$.  Moreover $ -\frac{1}{h_n}(v(z)-v(x))\ge f(x)> 0$, hence $m=v(x)>v(z)$ which contradicts the definition of $m$.
 \end{Proof}
 \begin{Proposition}\label{existenceD}
Let  $g:\G\to \R$ be such that
 \begin{equation}\label{compatD}
    g(x)\le \inf\left\{\sum_{k=0}^N h_n f(x_k)+g(y)\right\}\quad \forall x,y\in \G,
 \end{equation}
where the infimum is  taken over all the finite paths $\{x_0=x,x_1,\dots, x_N=y\}$ with $x_i\sim x_{i+1}$, $i=0,\dots,N-1$.
Then, the unique solution  of \eqref{HJD}  with the boundary condition
  \begin{equation}\label{BCD}
    u=g\qquad\text{ on $ \G$}
\end{equation}
is given by
 \begin{equation}\label{rappresentD}
  u(x)=\min\left\{ \sum_{k=0}^N h_nf(x_k)+g(y)\right\}
 \end{equation}
where the minimum is taken over all $y\in\G$ and over all the finite paths $\{x_0=x,x_1,\dots, x_N=y\}$ with $x_i\sim x_{i+1}$, $i=0,\dots,N-1$.
\end{Proposition}
 \begin{Proof} We observe that the function $u$ is well defined (i.e., the minimum is achieved) because, for $x$ fixed, the set of admissible paths is finite. Uniqueness is an immediate consequence of Prop.\ref{comparisonD}.
 In order to show that $u$ is a subsolution, given $x\in V_0^n$, for each $z\sim x$   consider a path $\{x_0=z,x_1,\dots, x_N=y_z\}$ such that
 $u(z)=g(y_z)+ \sum_{k=0}^N h_nf(x_k)$. Then $\{x, x_0,\dots, x_N\}$ is a   path connecting  $x$ to $\G$ and therefore
 \begin{align*}
    -(u(z)-u(x))&\le-\left(g(y_z)+h_n\sum_{k=0}^N f(x_k)-(g(y_z)+h_nf(x)+h_n\sum_{k=0}^N f(x_k))\right)\\
   &= h_nf(x).
 \end{align*}
 It follows that $|Du(x)|_n\le f(x)$.\par
 In order to prove that $u$ is a supersolution, given  $x\in V_0^n$  consider a path $\{x_0=x,x_1,\dots, x_N=y\}$ with $y\in \G$ such that $u(x)=g(y)+ \sum_{k=0}^N h_nf(x_k)$.  Then  $x_1\sim x$  and $\{x_1,\dots, x_N\}$ is a  path connecting  $x_1$ to $\G$. Therefore
 \begin{align*}
  -(u(x_1)-u(x))=-\left(u(x_1)-(g(y)+h_n\sum_{k=1}^N f(x_k))\right)+h_nf(x)\ge h_nf(x)
 \end{align*}
which implies  $|Du(x)|_n\ge f(x)$.\par
To show \eqref{BCD}, assume by contradiction that there exists $x\in \G$  and   a path $\{x_0=x,x_1,\dots, x_N=y\}$ with $y\in \G$
such that
\[u(x)=g(y)+h_n\sum_{k=0}^N f(x_k)<g(x).\]
By \eqref{compatD}, we have
\[g(y)+h_n\sum_{k=0}^N f(x_k)\ge g(y)+g(x)-g(y)=g(x)\]
and therefore a contradiction.
 \end{Proof}
\begin{Remark}\label{1:D} Observe that by \eqref{rappresentD}, the distance from the boundary
  \begin{equation}\label{distDB}
   d_n(x):=\inf\{d_n(x,y): \, y\in\G\} \qquad x\in V^n
 \end{equation}
is the unique solution of
\begin{equation*}
\left\{
\begin{array}{ll}
 |Du(x) |_n-1=0\quad& x\in V^n_0\\
 u=0\quad &  x\in\G.
 \end{array} \right.
\end{equation*}
\end{Remark}
\begin{Proposition}
Let $u$ be the solution of \eqref{HJD}-\eqref{BCD}. Then
\begin{align}
    &|u (x)|\le \max_{\G} |g|+  d_n(x)\,\max_{V^n}   \{f\}  &&\quad \forall x\in V^n\label{regD1}\\[6pt]
    &\frac{ |u(y)-u(x)|}{h_n}\le \max_{V^n}   \{f\} &&\quad \forall x,y\in V^n, x\sim y\label{regD2}
\end{align}
where $d_n(x)$ is the distance from the boundary introduced in \eqref{distDB}.
 \end{Proposition}
 \begin{Proof}
 The functions
\begin{equation*}
u^\pm(x):=\pm \left(\max_{\G} |g|+ d_n(x)\max_{V^n}  \{ f \}\right)
\end{equation*}
are a supersolution and respectively a subsolution of \eqref{HJD}-\eqref{BCD}, then the  estimate \eqref{regD1} follows  by Prop.\ref{comparisonD}.\par
 The estimate \eqref{regD2} is an immediate consequence of the equation \eqref{HJD}.
 \end{Proof}

Let $\bar u_n$ be a continuous piecewise linear reconstruction  of the solution $u_n$  of \eqref{HJD}-\eqref{BCD}  on $S^n$.
The previous proposition yields that the sequence $\{\bar u_{n}\}_n$ is compact. In Section \ref{HJSierpinski} we will show that it uniformly converges to a function $u$ and we will characterize $u$ in terms of an eikonal equation defined on the Sierpinski gasket.


\section{A continuous  eikonal equation on the prefractal}\label{s4}

In this section we introduce a continuous eikonal equation on the prefractal $S^n$. Moreover we estimate the distance
between the solutions  of the discrete problem \eqref{HJD} and of the corresponding continuous eikonal equations on $S^n$.
Let us recall that $V_0^n$ stands for the set of the interior vertices of $S^n$ while  $ \G$ is   the set of  the boundary vertices of $S^n$.
For any $n$,   fix   arbitrary orderings $I^n$  and $J^n$ of the vertices $x_i\in V^n$ and, respectively, of the edges $e_j\in E^n$ and
  denote by $\pi_j:[0,\ell_j]\to\R^2$,  $\ell_j>0$, a parametrization  of  the edge $e_j$.
For any $x_i \in V^n$, we set $Inc_{i}:=\{j\in J^n:x_i\in e_j\}$.
The   parametrization  $\pi_j$  of the arc  $e_j$ induces an orientation along  the edge, which can be expressed
by the \emph{signed incidence matrix} $A^n=\{a^n_{ij}\}_{i\in I^n,j\in J^n}$ with
\begin{equation}\label{incmatrix}
   a^n_{ij}:=\left\{
            \begin{array}{rl}
              1  & \quad\hbox{if $x_i\in e_j$ and $\pi_j(0)=x_i$,} \\
              -1 & \quad\hbox{if $x_i\in  e_j$ and $\pi_j(\ell_j)=x_i$,} \\
              0 & \quad\hbox{otherwise.}
            \end{array}
          \right.
\end{equation}
We denote by $ \d_n:S^n\times S^n\to \R^+$ the {\it path distance} on $S^n$, i.e.
 \begin{equation}\label{distNet}
   \d_n(x,y):=\inf\left\{\ell(\g):\,\text{$\g\subset S^n$ is a path joining $x$ to $y$ }\right\} \qquad \forall x,y\in S^n.
 \end{equation}
Note that $\d_n$ coincides with the vertex distance \eqref{geodistD} when restricted to $x,y\in V^n$. Moreover $\d_n\ge d$ and for any $n$ where $d$ as in \eqref{geodist}.\par
We say that $u$ is continuous in $S^n$ and write $u\in C(S^n)$ if $u$ is continuous with respect to the subspace topology of $S^n$.
In a similar way we define the space of upper semicontinuous functions  $\text{USC}(S^n)$ and  the space of lower semicontinuous functions $\text{LSC}(S^n)$, respectively.
For any function $u:S^n\to\R$, we denote by $u_j$ the restriction of $u$ to $e_j$, i.e.
\[u_j:=u\circ \pi_j:[0,\ell_j]\to \R.\]
Differentiation is defined with respect to the parametrization of the edge. Hence, if $x\in e_j$ for some $j\in J^n$ (i.e., $x$ is not a ramification point), we set
\[
D_{j} u(x):=\frac{du_j}{ d y} (y),\qquad\text{with $y=(\pi_j)^{-1}(x)$}
\]
while if $x_i\in V^n_0$ (i.e., $x$ is a vertex), we set
\[
D_{j} u(x_i):=\frac{du_j}{d y} (y)\quad\text{for any $e_j\in Inc_{i}$, $y=(\pi_j)^{-1}(x_i)$}.
\]
For a function $u:S^n\to\R$,  we consider the  eikonal equation
\begin{equation}\label{HJC}
      |Du| =f(x)   \qquad \textrm{in }  S^n.
\end{equation}
In the next definitions we introduce the class of admissible test functions and solution of \eqref{HJC} (see \cite{sc}).
\begin{Definition}\label{1:def2}
Let $\phi\in C(\G)$.
\begin{itemize}
  \item[i)] Let $x\in e_j$, $j\in J^n$. We say that $\phi$ is test function at $x$, if $\phi_j$ is differentiable at $\pi_j^{-1}(x)$.
  \item[ii)] Let $ x_i\in V^n_0$, $j,k\in Inc_i$, $j\neq k$. We say that $\phi$ is $(j,k)$-test function at $x$, if
  $\phi_j$ and $\phi_k$ are differentiable at $ \pi_j^{-1}(x)$ and $ \pi_k^{-1}(x)$, with
\begin{equation*}
a^n_{ij}D_j \phi(\pi_j^{-1}(x)  )+a^n_{ik}D_k \phi(\pi_k^{-1}(x) )=0,
\end{equation*}
where $(a^n_{ij})$ as in \eqref{incmatrix}.
\end{itemize}
\end{Definition}

\begin{Definition}\label{scsol}
\begin{itemize}
 \item[i)] If  $x\in e_j$, $j\in J^n$, then
a function $u\in\text{USC}(S^n)$ (resp., $v\in\text{LSC}(S^n)$) is called a   subsolution (resp. supersolution) of \eqref{HJC} at $x$ if  for any test function  $\phi$     for which  $u-\phi$ attains a local maximum at $x$ (resp., a local minimum), we have
\[
|D_j\phi(x)|\le f(x)\qquad    \big(\text{resp.,}\quad |D_j\phi(x)|\ge f(x)\big);
\]
\item [ii)] If $x=x_i\in V^n_0$, then
\begin{itemize}
\item[$\bullet$] A function $u\in\text{USC}(S^n)$ is called a   subsolution at $x$ if
for any $j,k\in Inc_i$ and any $(j,k)$-test function  $\phi$     for which
  $u-\phi$ attains a local maximum at $x$ relatively to $ e_j\cup   e_k$, we have
 $|D_j\phi(x)|\le f(x).$
  \item[$\bullet$] A function $v\in\text{LSC}(S^n)$ is called a  supersolution  at $x$ if  for any  $j\in Inc_i$, there exists $k\in Inc_i\setminus\{j\}$ (said feasible for $j$ at $x$) such that for any  $(j,k)$-test function  $\phi$
   for which $u-\phi$ attains a local minimum at $x$ relatively to $  e_j\cup   e_k$, we have
 $ |D_j\phi(x)|\ge f(x).$
 \end{itemize}
\end{itemize}
\end{Definition}
\begin{Remark}
It is worth to observe that the definitions of sub- and supersolution are not symmetric at the vertices. As observed in \cite{sc} for the equation $|D u|=1$, a definition of supersolution similar to the one of subsolution would not characterize the ``expected'' solution, i.e. the distance from the boundary. Also in \cite{ghn}, these definitions are asymmetric: see Def. \ref{gigasol} below.
Furthermore, let us recall that several definitions of solutions for Hamilton-Jacobi equations on network have been introduced recently; however for equation \eqref{HJC} they coincide (see \cite{cm_comp13}).

\end{Remark}
\begin{Remark}\label{standardef}
Consider a vertex $v_i\in V^n_0$ and, wlog, assume $\pi_j(0)=v_i$ $\forall j\in Inc_i$. Then, a function $u\in\text{USC}(S^n)$ is a viscosity subsolution in $v_i$ if, and only if, for every $j,k\in Inc_i$ with $j\neq k$, the function
\begin{equation*}
U_{jk}(t):=\left\{\begin{array}{ll}
u(\pi_j(t))&\qquad t\in[0,l_j)\\ u(\pi_k(-t))&\qquad t\in[-l_k,0)
\end{array}\right.\end{equation*}
is a (standard) viscosity subsolution at $t=0$.

Similarly, a function $u\in\text{LSC}(S^n)$ is a viscosity subsolution in $v_i$ if, and only if, for every $j\in Inc_i$ there exists $k\in Inc_i\setminus\{j\}$ such that the function $U_{jk}$ is a (standard) viscosity supersolution at $t=0$.
\end{Remark}
In the next propositions we recall from \cite{sc}   comparison  and existence result  for \eqref{HJC}; for the proof we refer the reader to \cite[Prop.6.1]{sc}.
 \begin{Proposition}\label{comparisonC}
Let $u$ and $v$ be  a subsolution and, respectively,  a supersolution of \eqref{HJC}  such that $u\le v$ on $\G$. Then $u\le v$ in $S^n$.
 \end{Proposition}
\begin{Proposition}\label{rapprC}
Let  $g:\G\to \R$ be such  that
 \begin{equation}\label{compatC}
    g(x)\le \inf\left\{\int_0^T f(\xi(t))dt+g(y)\right\}\quad \forall x,y\in \G
 \end{equation}
where the infimum is  taken over all the piecewise differentiable paths $\xi$ such that $\xi(0)=x$, $\xi(T)=y$.
 Then the unique solution  of \eqref{HJC}  with the boundary condition
  \begin{equation}\label{BCC}
    u=g\qquad\text{ on $ \G$}
\end{equation}
is given by
\begin{equation}\label{rappresentNet}
  u(x)=\inf\left\{\int_0^T f(\xi(t))dt+g(y)\right\}\quad\text{for all $x\in S^n$},
\end{equation}
where
the infimum is taken over all $y\in\G$ and over all the  paths $\xi$ such that $\xi(0)=x$, $\xi(T)=y$. Moreover, $u$ is bounded and Lipschitz continuous and
\begin{align}
   & |u (x)|\le \max_{\G} |g|+  \max_{S^n}  |f| \d_n(x,\G)\quad\text{for all $x\in S^n$}\label{regC1}\\
    &\|D u\|_\infty \le   \max_{S^n}  |f |. \label{regC2}
\end{align}
\end{Proposition}
\begin{Remark}\label{1:C}
Observe that by \eqref{rappresentNet}, the distance from the boundary
\begin{equation}\label{3.6c}
\d_n(x):=\inf\{\d_n(x,y): \, y\in \G\} \qquad x\in S^n
\end{equation}
is the unique solution of
\begin{equation*}
|Du(x)|=1\quad x\in V^n_0,\qquad  u=0\quad   x\in\G.
\end{equation*}
\end{Remark}
\begin{Remark}\label{rmk:compDC}
Note that assumption \eqref{compatD} for every $n\in\N$ ensures assumption \eqref{compatC}.
\end{Remark}
\begin{Proposition}\label{convC}
The sequence $\{ u_{n}\}_n$ of the solution of \eqref{HJC}-\eqref{BCC}   is decreasing and converges uniformly to a  continuous function $\bar u$  on $S$.
\end{Proposition}
\begin{Proof}
Let  $u_m$ be  the solution of \eqref{HJC}-\eqref{BCC} on the set $S^m$, with $m>n$. Since $S^n\subset S^m$,     $u_{m}$ is a subsolution of  \eqref{HJC}-\eqref{BCC} on $S^n$. By Proposition \ref{comparisonD}, it follows that $u_{m}\le u_{n}$. By \eqref{regC1}-\eqref{regC2}, the sequence $\{  u_{n}\}_n$ is equi-bounded, equi-continuous  and decreasing. Therefore  it converges  uniformly   to a function $\bar u$ on $S$.
\end{Proof}
\begin{Proposition}\label{stimaDC}
Let  $u_{h_n}$ and $u_n$ be respectively the solution of \eqref{HJD}-\eqref{BCD}  and of \eqref{HJC}-\eqref{BCC}. Then, there holds
\begin{equation*}
|u_n(x)-u_{h_n}(x)|\le C \og_f(h_n^{1/2})\qquad \forall x\in V^n
\end{equation*}
where $\og_f$ is the modulus of continuity of $f$.
Hence the sequences $u_{h_n}$ and $u_n$ converge to the same limit $\bar u$ on $S$.
\end{Proposition}
In order to prove this result, it is expedient to establish a preliminary lemma whose proof is postponed at the end of this section.
\begin{Lemma}\label{LemmaHJDK}
Let  $u_{h_n}$ be  the solution of \eqref{HJD}-\eqref{BCD}. Then $w_{h_n}=1-e^{-u_{h_n}}$
solves
\begin{equation}\label{HJDK}
\left\{
\begin{array}{ll}
|Dw|_n+ f_n(x)w=f_n(x)\quad &\textrm{in } V^n_0\\[6pt]
w=1-e^{-g}& \textrm{on } \G
\end{array}
\right.
\end{equation}
where $f_n(x)=h_n^{-1}( e^{h_nf(x)}-1)$. Moreover
\begin{equation}\label{propf1}
    \max_{S^n}|f_n-f|\le h_n \max_S  f.
\end{equation}
\end{Lemma}
\begin{Proofc}{Proof of Proposition \ref{stimaDC}}
We observe that the function $w_n:=1-e^{-u_n}$ solves
\[
|Dw| + f w=f\quad \textrm{in } S^n,\qquad
w=1-e^{-g}\quad \textrm{on } \G.
\]
Define $\Psi:V^n\times S^n \to \R$ by
\[\Psi(x,y)=w_{h_n}(x)-w_n(y)-\frac{\d_n(x,y)^2}{2\e}\]
where $w_{h_n}$ is the function introduced in Lemma~\ref{LemmaHJDK} while $\d_n$ is the geodesic distance on $S^n$ defined in \eqref{distNet}. Since $\Psi$ is continuous and $V^n\times S^n$ is compact, there exists
$(x_\e,y_\e)\in V^n\times S^n$ such that $\Psi(x_\e,y_\e)\ge\Psi(x,y)$ for any $(x,y)\in V^n\times S^n$. By $\Psi(x_\e,y_\e)\ge\Psi(x,x)$, we get
\begin{align*}
    \frac{\d_n(x_\e,y_\e)^2}{2\e}\le 2 (\|w_{h_n}\|_\infty+\|w_n\|_\infty)
\end{align*}
and therefore $\lim_{\e \to 0^+}\d_n(x_\e,y_\e)=0$. Hence there exists $x_i\in V^n$ such that
\begin{equation*}
   \lim_{\e\to 0^+} x_\e= \lim_{\e\to 0^+} y_\e=x_i.
\end{equation*}
Moreover, since $V^n$ is a finite set, wlog we can assume that, for $\e$ sufficiently small $x_\e=x_i$. By the Lipschitz continuity of $u_n$  and $\Psi(x_\e,y_\e)\ge\Psi(x_\e,x_\e)$, it follows
\[\frac{\d_n(x_\e,y_\e)^2}{2\e}\le w_n(x_\e)-w_n(y_\e) \le L_{w_n} \d_n(x_\e,y_\e)\]
where $L_{w_n}$ is the Lipschitz constant of $w_n$ (see Proposition \ref{rapprC}). Whence, we have
\begin{equation}\label{est2}
\d_n(x_\e,y_\e)\le L_{w_n}\e.
\end{equation}
Assume first that $x_\e=x_i\in \G$. Then, since $w_n=w_{h_n}$ on $\G$, we get
\begin{equation}\label{est3}
   w_{h_n}(x_\e)- w_n(y_\e)=w_n(x_\e)- w_n(y_\e)\le L_{w_n} \d_n(x_\e,y_\e).
\end{equation}
Assume $x_\e=x_i \in   V^n_0$ and let $z_\e\in V^n$ be such that $y_\e$ is contained in the edge $ e_j$, $j\in Inc_{i}$, of endpoints $x_\e$ and $z_\e$.
Since $z_\e\sim x_\e$  and  $\Psi(x_\e,y_\e)\ge\Psi(z_\e,y_\e)$, it follows
\begin{equation}\label{est4}
 -   \frac{1}{h_n}\left(\frac{\d_n(z_\e,y_\e)^2}{2\e}-\frac{\d_n(x_\e,y_\e)^2}{2\e}\right)+f_n(x_\e)w_{h_n}(x_\e)\le f_n(x_\e).
\end{equation}
By  $\d_n(x_\e, y_\e)+\d_n(y_\e, z_\e)=h_n$, we get
\begin{align*}
  -\left( \frac{\d_n(z_\e,y_\e)^2}{2\e}-\frac{\d_n(x_\e,y_\e)^2}{2\e}\right)= h_n \frac{\d_n(x_\e,y_\e)}{\e}-\frac{h_n^2}{ 2\e}
\end{align*}
and by \eqref{est4}
\begin{equation}\label{est5}
  \frac{\d_n(x_\e,y_\e)}{\e}  - \frac{h_n}{2\e}+f_n(x_\e)w_{h_n}(x_\e)\le f_n(x_\e).
\end{equation}
Defined $\phi(y)= w_{h_n} (x_\e)-\d_n(x_\e,y)^2/2\e$, by $-\Psi( x_\e,y_\e)\le -\Psi(x_\e,y)$  it follows that  $w_n-\phi$ attains a minimum at $y_\e$.
Observe that only two cases may occurr: either $y_\e=x_\e$ or  $y_\e\in e_j$. Let us consider the former case; the latter one can be dealt with by easy adaptations and we shall omit it.
Since $D_j\phi(x_\e)=0=D_k\phi(x_\e)$, we obtain that for any $k\in Inc_{x_\e}\setminus\{j\}$, $\phi$ is a $(j,k)$-admissible test  function for $w_n$ at $x_\e$.
By definition of viscosity supersolution, we obtain
\begin{equation}\label{est6}
    \frac{1}{\e}d(x_\e,y_\e)+f(y_\e)w_n(y_\e)\ge f(y_\e).
\end{equation}
By \eqref{est5} and \eqref{est6} we get
\begin{equation}\label{est7}
f_n(x_\e)w_{h_n}(x_\e)-    f(y_\e)w_n(y_\e)\le f_n(x_\e)-f(y_\e)+\frac{h_n}{2\e}.
\end{equation}
By either \eqref{est3} or \eqref{est7}, using \eqref{est2}, \eqref{regD1} and \eqref{propf1} (and our assumptions as well), we get
\begin{equation}\label{est8}
   w_{h_n}(x_\e)- w_n(y_\e)\le \frac{C}{\l}\left(\og_f(\e)+\frac{h_n}{2\e}\right)
\end{equation}
where $C$ depends only on $\max_\G|g|$,  $\max_{x\in V^n} d(x,\G)$, $\max_S f$  (note that $L_{w_n}$ depends only on $\max_S f$).
Taking $\e=h_n^{1/2}$ in \eqref{est8} we get the estimate
\[
w_{h_n}(x_\e)-w_{n}(y_\e)\le C(\og_f( h_n^{1/2})+h_n^{1/2})\le   C \og_f( h_n^{1/2}).
\]
Finally the inequality $\Psi(x_\e,y_\e)\ge\Psi(x,x)$ for $x\in V^n$ yields
\[u_{h_n}(x)-u_{n}(x)\le u_{h_n}(x_\e)-u_{n}(y_\e)\le C \og_f( h_n^{1/2}).\]
The proof of inequality  $u_{n}(x)-u_{h_n}(x)\le C \og_f( h_n^{1/2})$ can be obtained in a similar way; hence, we shall omit it.
\end{Proofc}
\begin{Proofc}{Proof of Lemma \ref{LemmaHJDK}}
For $x\in V^n_0$, consider $z\sim x$ such that
\[|Du_{h_n}(x)|_n=h_n^{-1}\left(u_{h_n}(z)-u_{h_n}(x)\right)=f(x).\]
Whence, by the definition of $w_{h_n}$, we deduce
\[\log\left(\frac{1-w_{h_n}(z)}{1-w_{h_n}(x)}\right)=h_nf(x)\]
and, by easy calculations, we get
\[-\left(w_{h_n}(z)-w_{h_n}(x)\right)-w_{h_n}(x)\left[1-e^{h_nf(x)}\right]=e^{h_nf(x)}-1.\]
Dividing by $h_n$ and taking into account the definition of $|Dw_{h_n}|_n$, we accomplish the proof.
\end{Proofc}
\section{The eikonal equation on the Sierpinski gasket}\label{HJSierpinski}

In this section we show that, in analogy with the construction  of the   Laplacian on fractal sets, the   eikonal equation on the Sierpinski gasket can be obtained as the limit of a sequence  of discrete problems defined on the prefractal $S^n$.\par
We first recall some notations and definitions in \cite{ghn}. Let $(\cX,d)$ be a metric space. Given a curve $\xi:I\subset\R\to \cX$, define the metric derivative of $\xi$ by
\[|\xi'|(t):=\lim_{s\to t}\frac{d(\xi(s),\xi(t))}{|s-t|}.\]
A curve $\xi$ is said {\it absolutely continuous} if $|\xi'|(t)$ exists for a.e. $t\in I$, $|\xi'|\in L^1_{loc}(I)$ and
\[d(\xi(s),\xi(t))\le \int_s^t|\xi'|(r)dr\]
where $d$ is the geodesic distance on $S$ defined in \eqref{geodist}.
For $\O\subset\cX$, denote by   $\cA(I,\O)$   the set of absolutely continuous curves such that $\xi(t)\in \O$ and $|\xi'|(t)\le 1 $ a.e. in $I$ and by $\cA_x(I,\O)$
the set of curve in  $\cA(I,\O)$ such that $0\in I$ and $\xi(0)=x$. For a   $\xi\in  \cA_x(I,\O)$ define the
exit and the entrance time in $\O$ by
\begin{align*}
    &T^{+}_\O[\xi]:=\inf\{t\in [0,+\infty):\,\xi(t)\in \pd \O\}\in [0,+\infty]\\
    &T^{-}_\O[\xi]:=\sup\{t\in (-\infty,0]:\,\xi(t)\in \pd \O\}\in [-\infty,0].
\end{align*}
A function $u$ is said arcwise upper (resp., lower) semicontinuous in $\O\subset\cX$ if for each $\xi\in \cA(I,\O)$
the function $u\circ \xi$ is upper (lower) semicontinuous in $I$. The set of arcwise upper (resp., lower) semicontinuous functions in $\O$ is denoted by $\text{USC}_a(\O)$ (resp., $\text{LSC}_a(\O)$). It is worth to recall that these definitions of semicontinuity are weaker than the corresponding standard ones (which are given in terms of the distance; see \cite{ghn}). For a function $w:\R\to\R$,  denote by $D^{\pm} w(t)$ respectively its super- and subdifferential at the point $t$ (see \cite{bcd} for the precise definition and main properties).

For $\O \subset\cX$, consider the eikonal equation
\begin{equation}\label{HJ}
    |Du|=f(x)\qquad \textrm{in } \O;
\end{equation}
for the sake of completeness, we recall the definition of (metric viscosity) subsolutions and supersolutions introduced in \cite{ghn}.
\begin{Definition}\label{gigasol}\hfill
\begin{itemize}
\item[i)] A function $u\in \text{USC}_a(\O)$ is said a viscosity subsolution of \eqref{HJ}
if for each $x\in \O$ and for all $\xi\in \cA_x(I,\O)$, the function $w:I\rightarrow \R$, with $w(t):=u(\xi(t))$, satisfies
\[|p|\le f(x)\qquad\forall p\in D^+w(0). \]
\item[ii)]
A function  $v\in \text{LSC}_a(\O)$ is said a viscosity supersolution of \eqref{HJ}
if for each $x\in \O$ and  $\e>0$, there exists $\xi\in \cA_x(\R,\O)$ and $w\in \text{LSC}(\R)$  such that
\[
\left\{
\begin{array}{l}
T^\pm:=T^{\pm}_\O[\xi]\textrm{ are both finite},\\
w(0)=v(x),\quad
w(t)\ge v(\xi(t))-\e\quad\forall t\in(T^-,T^+),\\
|p|\ge f(x)-\e \qquad\forall p\in D^-w(t),\, t\in(T^-,T^+).
\end{array}
\right.
\]
\end{itemize}
\end{Definition}
\begin{Remark}
The previous definitions apply to the case where $\cX$ is the Sierpinski gasket $S$ introduced in section \ref{sect:SG} with $\O:=S\setminus \G$ (recall that $\G$ is the set of the vertices of the starting simplex). In this case the metric on $\cX$ is the one induced by  the geodetic distance $d$ defined in \eqref{geodist}. Moreover arcwise  upper and lower semicontinuous functions coincide with upper and lower semicontinuous functions with respect to the metric induced by $d$.
\end{Remark}

The following results are an immediate consequence of the theory developed   in \cite{ghn} for a generic metric space $(\cX,d)$; for the proofs, we refer the reader to \cite[Thm3.1]{ghn} and to \cite[Thm3.2]{ghn}.
 \begin{Proposition}\label{comparisonCS}
Let $u$ and $v$ be  a subsolution and, respectively,  a supersolution of \eqref{HJ}  such that $u\le v$ on $\G$. Then $u\le v$ in $S$.
 \end{Proposition}
\begin{Proposition}
Let  $g:\G\to \R$ be such  that
 \begin{equation}\label{compatS}
    g(x)\le \inf\left\{\int_0^T f(\xi(t))dt+g(y)\right\}\quad\text{for all $x,y\in \G$}
 \end{equation}
where the infimum is  taken over all the paths $\xi\in\cA((0,T),S)$ such that $\xi(0)=x$, $\xi(T)=y$.
 Then the unique solution  of \eqref{HJ}  with the boundary condition
  \begin{equation}\label{BCS}
    u=g\qquad\text{ on $ \G$}
\end{equation}
is given by
 \begin{equation}\label{rappresentS}
  u(x)=\inf\left\{\int_0^T f(\xi(t))dt+g(y)\right\}\quad\text{for all $x\in S$}
 \end{equation}
where
the infimum is taken over all $y\in\G$ and over all  the paths $\xi\in\cA((0,T),S)$ such that $\xi(0)=x$, $\xi(T)=y$. Moreover
$u$ is continuous,  bounded  and
\begin{align*}
     |u (x)|\le \max_{\G} |g|+  \max_{S}  |f| d(x,\G)\quad\text{for all $x\in S^n$}.
\end{align*}
\end{Proposition}

The  Definition \ref{gigasol} also applies to the prefractal $S^n$.
In the next propositions we  establish    the relation between sub and supersolutions in the sense of Definition~\ref{gigasol} and of Definition~\ref{scsol} for the problem
\begin{equation}\label{gigapre}
|Du|=f(x)\qquad \forall x\in \O:=S^n\setminus \G.
\end{equation}
 These results are expedient in order to show the convergence of the discrete problem.

\begin{Proposition}
A function $u\in \text{USC}(S^n)$ is a subsolution in the sense of Definition \ref{gigasol} to problem \eqref{gigapre} if  and only  if it is a subsolution in the
sense of Definition \ref{scsol} to \eqref{gigapre}.
%
\end{Proposition}

Taking advantage of Remark \ref{standardef}, the proof can be easily done adapting to the case of $S^n$ the argument in \cite[Prop 6.1]{ghn} which shows the equivalence between metric viscosity subsolutions and standard viscosity subsolution in the Euclidean space. Therefore, we shall omit it.\par

\begin{Proposition} For $v\in \text{LSC}(S^n)$, the following statement are equivalent:
\begin{itemize}
\item[($i$)] $v$ is a supersolution in the sense of Definition \ref{gigasol} to problem \eqref{gigapre},
\item[($ii$)] $v$ is a supersolution in the sense of Definition \ref{scsol} to problem \eqref{gigapre},
\item[($iii$)] for each $x\in S^n$, there exists $\xi\in \cA_x(\R,S^n)$, with $T^\pm:=T^{\pm}_{\O}[\xi]$ both finite, such that the function $w(t):= v(\xi(t))$ $\forall t\in(T^-,T^+)$
    \[
|p|\ge f(\xi(t))  \qquad\forall p\in D^-w(t),\, t\in(T^-,T^+).
\]
\end{itemize}
\end{Proposition}
\begin{Remark}\label{rmkovvio}
Point ($iii$) coincides with the notion of supersolution in Definition~\ref{gigasol} with $\e=0$ and with $w:=w\circ \xi$.
\end{Remark}
\begin{Proof}
\emph{($i$) implies ($ii$).}
Assume by contradiction that $v$ does not satisfy the definition of supersolution in Def. \ref{scsol} at $x\in S^n $. We assume that $x$  coincides vith a vertex $x_i$   otherwise we can proceed as in \cite[Prop 6.2]{ghn}. Then there exists $j\in Inc_i$ such that for any $k\in Inc_i$, $k\neq j$, there exists a
$(j,k)$-admissible test function  $\phi_k$ for $u$ at $x$ which satisfies $|D_j\phi_k(x)|\le f(x)-2\d$ for some $\d>0$.
By adding the term $d(x,y)^2$ to the test function it is not restrictive to assume that
\[
    (v-\phi)(y)\ge d(x,y)^2\qquad \forall y\in S^n.
\]
By the regularity of $\phi_k$ we can find $r>0$ such that
\[
     |D \phi_k(y)|\le f(y)- \d \qquad \forall y\in B(x,2r)\cap(e_j\cup e_k).
\]
Define  $  \phi: \bigcup_{k\in Inc_i}  e_k\to \R$ by
\[  \phi(y):=\left\{
                  \begin{array}{ll}
                    \max_{k\in K} \phi_k(y), & \hbox{if $y\in  e_j$,} \\
                     \phi_k(y), & \hbox{if $y\in  e_k$, $k\in Inc_i$, $k\neq j$.}
                  \end{array}
                \right.
\]
As in \cite[Prop.4.1]{sc}, we have: $\phi$ is $C(\cup_{k\in Inc_i} e_k)$, it satisfies $\phi(x)=v(x)$ and (in viscosity and a.e. sense)
\[|D  \phi (y)|\le f(y)- \d \qquad \text{for all $y\in B(x,2r)$}.\]
For $\e\in (0,r^2/3)$, consider a couple $(\xi_\e,w_\e)$ as in the Definition \ref{gigasol} with $T^{\pm}=T^{\pm}_{B(x,r)}$ where $B(x,r)=\{y\in S^n:\,d(x,y)\le r\}$.
Let $\rho_\eta:\R\to\R$ be  a mollifier and set $\phi_{\eta,\e}(t)=( \phi (\xi_\e(\cdot))*\rho_\eta)(t)$.
Then, as $\eta\to 0$,  $\phi_{\eta,\e}(t) \to \phi(\xi_\e(\cdot))$ uniformly in $[T^-,T^+]$. Let $t_{\eta,\e}$ be a minimum point for $w_\e-\phi_{\eta,\e}$ in $[T^-,T^+]$. Then for $\eta$ sufficiently small in such a way that $\|\phi(\xi_\e(\cdot))-\phi_{\eta,\e}(\cdot)\|_{\infty,[T^-,T^+]}\le \e$, we have
\begin{eqnarray*}
0&=&v(x)-\phi(x) = w_\e(0)-\phi(\xi_\e(0))\ge w_\e(0)-\phi_{\eta,\e}(0)-\e\\
 &\ge& w_\e(t_{\eta,\e})-\phi_{\eta,\e}( t_{\eta,\e})-\e\ge
v (\xi_{ \e}(t_{\eta,\e}))-\phi (\xi_{ \e}(t_{\eta,\e}))-3\e\\&\ge&  d(\xi_\e(t_{\eta,\e}),x)^2-3\e.
\end{eqnarray*}
By our choice of $\e$, $\xi_\e$ belongs to $B(x,r)$ for $t\in[0,t_{\eta,\e}]$; hence, $t_{\eta,\e}$ belongs to $(T^-,T^+)$ and $\phi_{\eta,\e}$ is a test function for $w_\e$ at $t_{\eta,\e}$. It follows
that
\begin{eqnarray*}
f(\xi_\e(t_{\eta,\e}))-\e&\le& |D\phi_{\eta,\e}(t_{\eta,\e})|=|D(\phi(\xi_\e(\cdot)*\rho_\eta)(t_{\eta,\e})|\\
&\le&
\int_{(t_{\eta,\e}-\eta,t_{\eta,\e}+\eta)}|D \phi(\xi_\e(s))||\dot \xi_\e(s)|\rho_\eta(t_{\eta,\e}-s)ds\\
&\le&
\int_{(t_{\eta,\e}-\eta,t_{\eta,\e}+\eta)}f(\xi_\e(s))\rho_\eta(t_{\eta,\e}-s)ds-\d\\
&\le& f(\xi_\e(t_{\eta,\e}))+\og_f(\eta)-\d
\end{eqnarray*}
which gives a contradiction for $\e$ and $\eta$ small.\\

\emph{($ii$) implies ($iii$).} Let us now prove that if $v$ is a supersolution in the sense of Def. \ref{scsol}, then it is a supersolution in the sense of Def.  \ref{gigasol}. To this end, fix a point $x\in S^n$ and introduce a map $\xi:\R\to S^n$ as follows. Assume that $x\in e_j$ with $x=\pi_j(y_x)$, $y_x\in [0,l_j]$. Set $\xi(t):=\pi_j(y_x+t)$ for $t\in[-y_x,l_j-y_x]$; in particular we have $\xi(0)=x$. Let us now proceed to introduce $\xi$ for positive $t$; the case of negative $t$ is similar and we shall omit it. The point $\xi(l_j-y_x)$ is an endpoint of $e_j$; hence it belongs either to $\G$ or to $V^n_0$. In the former case we obtain $T^+$. In the latter case, by the definition of supersolution  in Def \ref{scsol} there exists a feasible edge for $e_j$, say $e_k$; wlog assume $\pi_k(0)=\xi(l_j-y_x)$. We define $\xi(t):=\pi_k(t-l_j+y_x)$ for $t\in(l_j-y_x, l_k+l_j-y_x)$. We iterate this idea of following an edge and then choosing the corresponding feasible edge when arriving at a vertex of $V^n_0$.
Assume now that the initial point $x$ belongs to $V^n_0$, say $x=v_i$. In this case we choose arbitrarily an edge $e_j$, incident to $v_i$ and the corresponding feasible edge in $v_i$, say $e_k$. Wlog we assume $x=\pi_j(0)=\pi_k(0)$. We introduce $\xi(t)=\pi_j(t)$ for $t\in[0,l_j]$ and $\xi(t)=\pi_k(-t)$ for $t\in[-l_k,0)$ and we iterate the same idea as before. We claim that
\begin{equation}\label{claim:Tpm}
 \textrm{ $T^\pm:=T^\pm_{S^n}[\xi]$ are both finite.}
\end{equation}
In order to prove this property, we proceed by contradiction assuming that $\xi(t)\notin\G$ for every $t\in\R$. Since $S^n$ is given by a finite collection of edges of finite length, this may happen only if $\xi$ contains a loop, say $L:=  e_1\cup\dots\cup  e_k$.
Being lower semicontinuous, the function $v$ attains its minimum with respect to $L$ at some point $\bar x$. Since the constant function $\phi=\min_Lv$ is an admissible test function for $v$ at $\bar x$ we obtain a contradiction with the definition of supersolution in Def. \ref{scsol}; hence our claim \eqref{claim:Tpm} is completely proved.

Finally, by the Definition \ref{scsol} and Remark \ref{standardef}, one can easy check that for $w(t):=v(\xi(t))$ $\forall t\in [T^-,T^+]$, there holds: $|p|\ge f(\xi(t))$ for every $p\in D^-w(t)$, $t\in(T^-,T^+)$.

\emph{($iii$) implies ($i$).} By Remark \ref{rmkovvio}, this statement   is a straightforward consequence of the definition of supersolution by in Def.\ref{gigasol}. Hence, the proof is accomplished.
\end{Proof}

A natural question in the previous approach is if we need to consider  all the paths on the Sierpinski gasket or just the ones  that visit a finite number
of vertices. We say that a path $\xi:I\to S$ is  {\it piecewise differentiable} when it is absolutely continuous and there exist $n\in\N$, $t_0<t_1<\dots<t_m$ and $j_1,\dots,j_m\in J^n$ such that there holds: $I=(t_0,t_m)$ and $\xi(t)$ belongs to $e_{j_k}$ for $t\in (t_{k-1}, t_k)$. In other words, a piecewise differentiable path is an absolutely continuous path in some $S^n$ and it changes edge only a finite number of times.
\begin{Definition}
We say that a function $u$ is a {\it mild} subsolution (resp., supersolution) of \eqref{HJ} if it verifies the definition  of subsolution (resp., supersolution) in Def. \ref{gigasol} by  considering only piecewise differentiable paths.
\end{Definition}
\begin{Proposition}\label{ex_mild}
Assume \eqref{compatS}.  Then, the  unique mild  solution to \eqref{HJ}-\eqref{BCS} is given by
 \[
  v(x)=\inf\left\{\int_0^T f(\xi(t))dt+g(y)\right\}\quad\text{for all $x\in S$}
 \]
where the infimum is taken over all $y\in\G$ and over all the piecewise differentiable paths $\xi\in\cA((0,T),S)$ such that $\xi(0)=x$, $\xi(T)=y$.
Moreover $v$ coincides with the solution $u$ given by \eqref{rappresentS}.
\end{Proposition}
\begin{Proof}
Following exactly the same arguments in the proof of \cite[Thm3.1]{ghn}, one can easily check that Prop. \ref{comparisonCS} still holds for mild super- and subsolutions. Similarly, repeating the arguments of the proof of \cite[Thm3.2]{ghn}, one can also obtain that the function $v$ is a mild supersolution.\par
Finally, we observe that a subsolution (respectively, a supersolution) in the sense of Def.\ref{gigasol} is also a mild subsolution (resp.,  mild supersolution). Hence  we deduce that the function $u$ defined in \eqref{rappresentS} is a mild solution and the comparison principle yields that $u=v$.
\end{Proof}

\begin{Remark}\label{1:S}
Observe that by \eqref{rappresentS}, the distance from the boundary 
\begin{equation}\label{3.6s}
d(x):=\inf\{d(x,y): \, y\in \G\} \qquad x\in S
\end{equation}
is the unique solution of
\begin{equation*}
|Du(x) |=1\quad x\in S,\qquad  u=0\quad   x\in\G.
\end{equation*}
\end{Remark}
\begin{Theorem}\label{converg}
Assume \eqref{compatS} and
let $u_n$ be the viscosity solution of \eqref{HJC}-\eqref{BCC}.
Then, as $n\to+\infty$, $u_n$ tends to $u$ uniformly in $S$ where $u$ is the viscosity solution of \eqref{HJ}-\eqref{BCS}.
\end{Theorem}
\begin{Proof}
First observe that if condition \eqref{compatS} is satisfied, then, for any $n$, condition \eqref{compatC} is satisfied as well, hence for any $n$ there exists a solution to \eqref{HJC}-\eqref{BCC}.
By Prop. \ref{convC}, we already know that $u_n$ converges uniformly to a function $\bar u$ on $S$.\\
In order to show that $\bar u$ is a supersolution at $x\in S$, given $\e>0$, let $n\in \N$ be such that $x\in S^n$ and $0\le u_n(x)-u(x)\le \e$.  Let $\xi\in \cA_x(\R, S^n)$ and $w\in \text{LSC}(\R)$  be a $\e$-admissible couple for $u_n$ as in Definition  \ref{gigasol} with $\cX=S^n$. Set $\bar w=w+\bar u(x)-u_n(x)$. Then $\bar u(x)=\bar w(x)$ and $u(\xi(t))\le u_n(\xi(t))\le w(\xi(t))+2\e$. Moreover $|p|\ge f(\xi(t))-\e$ for all $p\in D^-  w(\xi(t))=D^- \bar w(\xi(t))$. It follows that $(\xi,w)$ is a  $\e$-admissible couple for $\bar u$ in the sense of Definition \ref{gigasol}.\par
To show that $\bar u$ is a subsolution, by Proposition \ref{ex_mild}, it suffices to show that it is a mild subsolution. To this end, given $x\in S$, consider  $\xi\in A_x(\R, S)$ and  $\bar n\in \N$ such that $\xi\in  \cA_x(\R,S^{\bar n})$. Set $w(t):=\bar u(\xi(t))$ and let $\phi\in C^1(\R)$ be supertangent to $w$ in $0$. Since $u_n \to u$ uniformly, then $w_n(t):=u_n(\xi(t))$ converges to  $w$ uniformly. Hence  there exists  $t_n \to 0$ such that   $\phi$ is supertangent to $w_n$ at $t_n$. Moreover $\xi\in \cA_x(\R,S^n)$ for $n>\bar n$ and  therefore $|D\phi(t_n)|\le f(\xi(t_n))$. For  $n\to \infty$,  since $\xi(t_n)\to\xi(0)$ it follows that $|D\phi(0)|\le f(x)$.
\end{Proof}

As an immediate consequence of  Prop. \ref{stimaDC}  and Thm. \ref{converg}, we get the convergence of discrete problems to the continuous problem on $S$.
\begin{Theorem}\label{convDS}
Assume \eqref{compatS} and \eqref{compatD} for any $n\in\N$.
Then, as $n\to+\infty$, the  sequence $\{u_{h_n}\}_n$  of the solutions of \eqref{HJD}-\eqref{BCD} converges uniformly to the solution $u$ of \eqref{HJ}-\eqref{BCS}.
 \end{Theorem}
 \begin{Remark}
Observe that \eqref{compatS} does not imply \eqref{compatD} except in some particular case (e.g., $g$ constant on $\G$).
Nevertheless we can always modify $f$ in such a way that \eqref{compatS} implies \eqref{compatD}. Indeed, under condition \eqref{compatS}, let us show that \eqref{compatD} is satisfied with $f$ replaced by $f_n(x):=f(x)+\og_f(h_n)$. For any $x,y\in\G$ and any path $\{x_0=x,x_1,\dots, x_N=y\}\subset V^n$ with $x_i\sim x_{i+1}$, $i=0,\dots,N-1$, define an admissible continuous path $\xi:[0, Nh_n]\to S^n$ such $\xi(ih_n)=x_i$, $i=0,\dots, N-1$. Then, by \eqref{compatS}, there holds
\begin{align*}
    g(x)\le \int_0^{Nh}f(\xi(t))dt +g(y)\le \sum_{i=0}^{N-1}\int_{ih}^{(i+1)h}f(\xi(t))dt +g(y)\\
        \le \sum_{i=0}^{N-1}h_n(f(x_i)+\og_f(h_n))+g(y)=\sum_{i=0}^{N-1}h_n f_n(x_i)+g(y).
\end{align*}
Therefore, 
$f_n$ satisfies \eqref{compatD}. Moreover, since $\{f_n\}_n$ is decreasing, by the same arguments as those of previous sections, one can prove that the solution $u_{h_n}$ to \eqref{HJD}-\eqref{BCD} with $f$ replaced by $f_n$ converges to the solution of \eqref{HJ}-\eqref{BCS}.
 \end{Remark}
\begin{Corollary}
Let $\d_n$, $d_n$ and $d$ be respectively the distance function from $\G$ respectively on $S^n$, $V^n$ and $S$, defined in \eqref{3.6c}, \eqref{distDB} and \eqref{3.6s}. Then, as $n\to+\infty$, both $\d_n$ and $d_n$ converge to $d$.
\end{Corollary}
\begin{Proof}
By Remark \ref{1:D} (respectively, Remarks \ref{1:C} and \ref{1:S}), the function $d_n$ (resp., $\d_n$ and $d$) is the solution of the eikonal equation with $f\equiv 1$ in $S^n$ (resp., $V^n$ and $S$) and $g=0$ on $\G$. Invoking Theorem \ref{converg} (resp., Thm.\ref{convDS}), we infer that $\d_n$ (resp., $d_n$) converge to $d$ as $n\to+\infty$.
\end{Proof}
%
\section{Extensions}\label{conclusion}
In this section we discuss some possible  extensions of the previous results. \par

\indent\textbf{More general Hamiltonians: }
For the sake of clarity, the present paper only concerns the eikonal equation; it is our purpose to extend these results to more general Hamiltonians in a future work.
For instance, we observe that the results on prefractals can be extended to an Hamiltonian $H(x,p)$ convex and coercive with respect to $p$,
obtaining the uniform convergence of the solutions of \eqref{HJD}-\eqref{BCD} to a function $u$ defined on the Sierpinski gasket.
The  definition  introduced in \cite{ghn} for the sole eikonal equation   can possibly be extended to a class of
Hamiltonians of the type $H(x,|p|)$.
\textbf{Boundary conditions: }
Following the classical definition of Laplacian on the Sierpinski gasket,
we imposed the boundary conditions on $\G$.
Nevertheless, by easy modifications, it is possible to consider as boundary set any finite subset of $V$. An interesting case considered in \cite{hs,s1} is when the boundary reduces to one of point of $\G$ and $g=0$ at this point. The solutions of the \eqref{HJD}, \eqref{HJC}  and \eqref{HJ} represent  respectively the vertex distance on $S^n$, the path distance on $S^n$  and the path distance on $S$ from the given boundary vertex.
More generally the problem \eqref{HJ}-\eqref{BCS} can be considered in a  connected subdomain $\O$ of $S$ imposing the boundary condition on  any finite subset of   vertices  contained in $\O$.\\

\textbf{More general fractals: }
The interior approximation method can be extended to the class of post-critically finite self similar sets (generally speaking, these sets are obtained by subdividing the initial cell into cells of smaller and smaller size and the cells must intersect at isolated points; see \cite{k1,s2} for the precise definition). We will consider the problem in a future work.\\

\indent\textbf{Resistance metric: }
Given the energy form associated to the Laplacian on the Sierpinski gasket, it is possible to introduce a metric, called resistance metric, which plays an important role in several estimates (see \cite{k1}, \cite{s2}). While in the Euclidean setting
resistance metric and pathwise distance coincide, this is not true on  $S$.  On the other hand the Sierpinski gasket can be embedded in $\R^2$ by certain harmonic maps whose image is now called the harmonic Sierpinski gasket and on this set the two  distances coincide (see \cite{ hi}, \cite{kz}). We aim to  consider  a characterization of the resistance metric via the  eikonal equation in a future work.

\section*{Acknowledgment}
This work is partially supported by European Union under the 7th Framework Programme FP7-PEOPLE-
2010-ITN Grant agreement number 264735-SADCO and by Indam-Gnampa.


\end{document}